\documentclass[11pt]{amsart}

\usepackage{amsmath}
\usepackage{amssymb}
\usepackage{graphicx}

\theoremstyle{definition}

\numberwithin{equation}{section}

\title[Solvability of quadratic integral systems]
{On the solvability of some systems of quadratic integral equations in dimensions two and three}

\author[V. Vougalter]{Vitali Vougalter}

\address[V. Vougalter]{Department of Mathematics,  University of Toronto,
Toronto, Ontario, M5S 2E4, Canada}
\email{{\tt vitali@math.toronto.edu}}



\begin{document}

\begin{abstract}
The work deals with the existence of solutions of a certain system
of quadratic integral equations in $H^{2}({\mathbb R}^{d}, {\mathbb R}^{N}), \ d=2,3$. 
We demonstrate the existence of a perturbed solution by virtue of
a fixed point technique.
\end{abstract}

\maketitle


\setcounter{equation}{0}

\section{Introduction}

The present article is devoted to the existence of 
solutions of the following system of integral equations,
\begin{equation}
\label{p}
u_{m}(x)=u_{0, m}(x)+[T_{m}u_{m}(x)]
\int_{{\mathbb R}^{d}}K_{m}(x-y)g_{m}(u(y))dy, 
\end{equation}
where $x\in {\mathbb R}^{d}, \ d=2,3,  \ 1\leq m\leq N$. The vector function contained in system (\ref{p}) is
given by
\begin{equation}
\label{vf}
u(x)=(u_{1}(x), u_{2}(x), ..., u_{N}(x))^{T}.  
\end{equation}
Analogously,
\begin{equation}
\label{vf0}
u_{0}(x)=(u_{0, 1}(x), u_{0, 2}(x), ..., u_{0, N}(x))^{T}.  
\end{equation}
The precise conditions on $u_{0}(x), \ g_{m}(u)$, the linear operators
$T_{m}$ and the kernels $K_{m}(x)$ will be specified further down. For each equation
of problem (\ref{p}) the second term in its
right side is a product of $T_{m}u_{m}(x)$ and the integral operator
acting on the function $g_{m}(u)$, for which the sublinear growth will be
established in the proof of Theorem 1.3. Therefore, the system of
integral equations of this kind is called {\it quadratic}. Existence of solutions of
the single integral equation of this kind was discussed in ~\cite{V24}.

The theory of 
integral equations has many important applications in describing various
events and phenomena of the real world. This is caused by the fact that this
theory is often applicable in different branches of mathematics, in mathematical
physics, economics, biology as well as in dealing with certain real world problems.
Quadratic integral equations appear in the theories of the radiative
transfer and neutron transport, in the kinetic theory of gases, in the design
of the bandlimited signals for the binary communication via the simple
memoryless correlation detection when signals are disturbed by 
additive white Gaussian noises (see e.g. ~\cite{AC08}, ~\cite{A92}, ~\cite{E05}
and references therein).
Work ~\cite{AC08} deals with the
solvability of a nonlinear quadratic integral equation in the Banach space
of real functions being defined and continuous on a bounded and closed
interval by means of the fixed point result. Articles ~\cite{AC11} and ~\cite{AC14}
are devoted to the existence of solutions for 
quadratic integral equations on unbounded intervals. The solvability
of quadratic integral inclusions was covered in ~\cite{AC13}.
In ~\cite{DH12} the authors investigate the nondecreasing solutions of a
quadratic integral equation of Urysohn-Stieltjes type.
The existence of solutions of
quadratic integral equations in Orlicz spaces was studied in
~\cite{CM12}, ~\cite{CM15}, ~\cite{CM16}. 
The integro-differential equations, which may contain either Fredholm or
non-Fredholm operators, appear in mathematical biology when
considering systems with nonlocal consumption of resources and 
intra-specific competition (see ~\cite{EV22}, ~\cite{VV21},
~\cite{VV22} and references therein). 
The contraction argument was applied
in ~\cite{V10} to estimate the perturbation to the standing solitary wave
of the Nonlinear Schr\"odinger equation when either the external
potential or the nonlinear term was perturbed. The similar ideas were exploited
to show the persistence of pulses for certain reaction-diffusion type
equations (see ~\cite{CV21}).

Suppose that the assumption below is fulfilled.

\medskip

\noindent
{\bf Assumption 1.1.}  {\it For $1\leq m\leq N$, the kernels
$K_{m}(x): {\mathbb R}^{d}\to {\mathbb R}, \ d=2,3$ are nontrivial, such that $K_{m}(x), \Delta K_{m}(x)\in L^{1}({\mathbb R}^{d})$.
The functions
$u_{0, m}(x): {\mathbb R}^{d}\to {\mathbb R}$ do not vanish identically in  ${\mathbb R}^{d}$ for some $1\leq m\leq N$ and 
$u_{0, m}(x)\in H^{2}({\mathbb R}^{d})$.
We also suppose that the linear operators
$T_{m}: H^{2}({\mathbb R}^{d})\to H^{2}({\mathbb R}^{d})$ are bounded, such that their
norms $0<\|T_{m}\|<\infty$.}  

\bigskip

Let us use the Sobolev space
\begin{equation}
\label{H1}  
H^{2}({\mathbb R}^{d}):=\big\{\phi(x):{\mathbb R}^{d}\to {\mathbb {\mathbb R}} \ | \
\phi(x)\in L^{2}({\mathbb R}^{d}), \ \Delta \phi(x)\in
L^{2}({\mathbb R}^{d}) \big \}
\end{equation}
with $d=2,3$.
It is equipped with the norm
\begin{equation}
\label{n}
\|\phi\|_{H^{2}({\mathbb R}^{d})}^{2}:=\|\phi\|_{L^{2}({\mathbb R}^{d})}^{2}+
\|\Delta \phi\|_{L^{2}({\mathbb R}^{d})}^{2}.
\end{equation}
For a vector function (\ref{vf}) we have the norm
\begin{equation}
\label{nvf}
\|u\|_{H^{2}({\mathbb R}^{d}, {\mathbb R}^{N})}^{2}:=\|u\|_{L^{2}({\mathbb R}^{d}, {\mathbb R}^{N})}^{2}+
\sum_{m=1}^{N}\|\Delta u_{m} \|_{L^{2}({\mathbb R}^{d})}^{2},
\end{equation}
where
\begin{equation}
\label{nvf2}
\|u\|_{L^{2}({\mathbb R}^{d}, {\mathbb R}^{N})}^{2}:=\sum_{m=1}^{N}\|u_{m}\|_{L^{2}({\mathbb R}^{d})}^{2}.
\end{equation}
Let us introduce another quantity relevant to our argument, namely
\begin{equation}
\label{11}
\|K_{m}\|_{{\tilde W}^{2, 1}({\mathbb R}^{d})}:=\sqrt{\|K_{m}\|_{L^{1}({\mathbb R}^{d})}^{2}+
\|\Delta K_{m}\|_{L^{1}({\mathbb R}^{d})}^{2}}
\end{equation}
and define
\begin{equation}
\label{Q}
Q:=\sqrt{\sum_{m=1}^{N}\|T_{m}\|^{2}\|K_{m}\|_{{\tilde W}^{2, 1}({\mathbb R}^{d})}^{2}},  
\end{equation}
so that $0<Q<\infty$.

It can be easily verified that the operators
\begin{equation}
\label{T}  
T_{m}\phi(x):=(-\Delta+1)^{-1}\phi(x), \quad \phi(x)\in H^{2}({\mathbb R}^{d})
\end{equation}
satisfy Assumption 1.1. 

Let us recall the Sobolev embedding
\begin{equation}
\label{e}
\|\phi(x)\|_{L^{\infty}({\mathbb R}^{d})}\leq c_{e}
\|\phi(x)\|_{H^{2}({\mathbb R}^{d})}, \ d=2,3,
\end{equation} 
where $c_{e}>0$ is a constant.
 
We will use the algebra property for our Sobolev space.
For any $\phi_{1}(x), \ \phi_{2}(x)\in H^{2}({\mathbb R}^{d}), \ d=2,3$, we have
$\phi_{1}(x)\phi_{2}(x)\in H^{2}({\mathbb R}^{d})$ and 
\begin{equation}
\label{alg}
\|\phi_{1}(x)\phi_{2}(x)\|_{H^{2}({\mathbb R}^{d})}\leq c_{a}\|\phi_{1}(x)\|_{H^{2}({\mathbb R}^{d})}
\|\phi_{2}(x)\|_{H^{2}({\mathbb R}^{d})},
\end{equation}
where $c_{a}>0$ is a constant. 

By virtue of the Young's inequality (see e.g.
Section 4.2 of ~\cite{LL97}), we have the upper bound on the norm of the
convolution as
\begin{equation}
\label{y}
\|\phi_{1}*\phi_{2}\|_{L^{2}({\mathbb R}^{d})}\leq \|\phi_{1}\|_{L^{1}({\mathbb R}^{d})}
\|\phi_{2}\|_{L^{2}({\mathbb R}^{d})}.
\end{equation}
Clearly, estimate (\ref{y}) yields the inequality
\begin{equation}
\label{yd}
\|\Delta_{x}\int_{{\mathbb R}^{d}}\phi_{1}(x-y)\phi_{2}(y)dy\|
_{L^{2}({\mathbb R}^{d})}\leq
\|\Delta \phi_{1}\|_{L^{1}({\mathbb R}^{d})}\|\phi_{2}\|_{L^{2}({\mathbb R}^{d})}.
\end{equation} 
We seek the solution of the nonlinear system 
(\ref{p}) of the form
\begin{equation}
\label{r}
u(x)=u_{0}(x)+u_{p}(x),
\end{equation}
where $u_{0}(x)$ is given by (\ref{vf0}) and
\begin{equation}
\label{up}
u_{p}(x):=(u_{p, 1}(x), u_{p, 2}(x), ..., u_{p, N}(x))^{T}.  
\end{equation}  
Evidently, we arrive at the perturbative system of equations
$$
u_{p, m}(x)=
$$
\begin{equation}
\label{pert}
[T_{m}(u_{0, m}(x)+u_{p, m}(x))]\int_{{\mathbb R}^{d}}
K_{m}(x-y)g_{m}(u_{0}(y)+u_{p}(y))dy, 
\end{equation}
with $1\leq m\leq N, \ d=2,3$.

We introduce a closed ball in the Sobolev space
\begin{equation}
\label{b}
B_{\rho}:=\{u(x)\in H^{2}({\mathbb R}^{d}, {\mathbb R}^{N}) \ | \
\|u\|_{H^{2}({\mathbb R}^{d}, {\mathbb R}^{N})}\leq \rho \}, \quad 0<\rho\leq 1.
\end{equation}
Let us look for the solutions of system (\ref{pert}) as fixed points of
the auxiliary nonlinear problem
\begin{equation}
\label{aux}
u_{m}(x)=[T_{m}(u_{0, m}(x)+v_{m}(x))] \int_{{\mathbb R}^{d}}
K_{m}(x-y)g_{m}(u_{0}(y)+v(y))dy,
\end{equation}
where $1\leq m\leq N, \ d=2,3$ in the ball (\ref{b}).

We define the closed ball in the space of $N$ dimensions as
\begin{equation}
\label{i}
I:=\Big\{z\in {\mathbb R}^{N} \ | \ |z|_{{\mathbb R}^{N}}\leq c_{e}
(\|u_{0}\|_{H^{2}({\mathbb R}^{d}, {\mathbb R}^{N})}+1)\Big\}
\end{equation}
along with the closed ball in the space of $C^{1}(I, {\mathbb R}^{N})$, such that
$$
D_{M}:=
$$
\begin{equation}
\label{M}
\{g(z):=(g_{1}(z), g_{2}(z), ..., g_{N}(z))\in C^{1}(I, {\mathbb R}^{N})
\ | \ \|g\|_{C^{1}(I, {\mathbb R}^{N})}\leq M \}
\end{equation}
with $M>0$.
In this context the norm
\begin{equation}
\label{gn1}
\|g\|_{C^{1}(I, {\mathbb R}^{N})}:=\sum_{m=1}^{N}\|g_{m}\|_{C^{1}(I)},
\end{equation}
\begin{equation}
\label{gn2}
\|g_{m}\|_{C^{1}(I)}:=\|g_{m}\|_{C(I)}+\sum_{n=1}^{N}\Big\|\frac{\partial g_{m}}
{\partial z_{n}}\Big\|_{C(I)},
\end{equation}  
where $\|\phi\|_{C(I)}:=\hbox{max}_{z\in I}|\phi(z)|$ for $\phi\in C(I)$.

\bigskip

\noindent
{\bf Assumption 1.2.} {\it Let $1\leq m\leq N$. It is assumed that
$g_{m}(z): {\mathbb R}^{N}\to {\mathbb R}$, such that
$g_{m}(0)=0$. We also assume that $g(z)|_{I}\in D_{M}$ and
it does not vanish identically in the ball $I$}.

\bigskip

Let $t_{g}$ be the operator defined by the right side of the system of equations (\ref{aux}),
such that $u = t_{g}v$.
Our first main result is as follows.

\bigskip

\noindent
{\bf Theorem 1.3.} {\it Let Assumptions 1.1 and 1.2 be valid and
\begin{equation}
\label{rub}
c_{a}M(\|u_{0}\|_{H^{2}({\mathbb R}^{d}, {\mathbb R}^{N})}+1)^{2}Q\leq \frac{\rho}{2}.
\end{equation}
Then the map $t_{g}: B_{\rho}\to B_{\rho}$ associated with problem
(\ref{aux}) is a strict contraction.
The unique fixed point $u_{p}(x)$ of this map $t_{g}$ is the only solution of
the system of equations (\ref{pert}) in $B_{\rho}$.}

\bigskip

Obviously, the resulting solution of system (\ref{p}) given by formula (\ref{r}) will
be nontrivial in ${\mathbb R}^{d}, \ d=2,3$ because all $g_{m}(0)=0$, the operators
$T_{m}$ are linear and the functions $u_{0, m}(x)$ do not vanish identically in the whole space for a certain
$1\leq m\leq N$ as we assume.
  
\bigskip

For the technical purpose we define the quantity
\begin{equation}
\label{sig}
\sigma:=2c_{a}QM(\|u_{0}\|_{H^{2}({\mathbb R}^{d}, {\mathbb R}^{N})}+1)>0. 
\end{equation}
Our second major statement deals with the continuity of the cumulative solution
of system (\ref{p}) given by (\ref{r}) with respect to the nonlinear
vector-function $g$.  

\bigskip

\noindent
{\bf Theorem 1.4.} {\it Suppose that the assumptions of Theorem 1.3 hold.
Let $u_{p,j}(x)$ be the unique fixed point of the map
$t_{g_{j}}: B_{\rho}\to B_{\rho}$ for $g_{j}, \ j=1,2$ and the resulting solution
of system (\ref{p}) with $g(z)=g_{j}(z)$ is given by
\begin{equation}
\label{u12}  
u_{j}(x)=u_{0}(x)+u_{p, j}(x).
\end{equation}
Then} 
\begin{equation}
\label{cont}
\|u_{1}(x)-u_{2}(x)\|_{H^{2}({\mathbb R}^{d}, {\mathbb R}^{N})}\leq
\end{equation}
$$
\frac{\sigma}
{2M (1-\sigma)}(\|u_{0}\|_{H^{2}({\mathbb R}^{d}, {\mathbb R}^{N})}+1)\|g_{1}(z)-g_{2}(z)\|_
{C^{1}(I, {\mathbb R}^{N})}.
$$

\bigskip

Let us proceed to the proof of our first main proposition.

\bigskip 


\setcounter{equation}{0}

\section{The existence of the perturbed solution}

\noindent
{\it Proof of Theorem 1.3.} We choose arbitrarily $v(x)\in B_{\rho}$.
By means of (\ref{aux}) along with (\ref{alg}), we obtain the upper bound
for $1\leq m\leq N, \ d=2,3$ as
\begin{equation}
\label{auxub} 
\|u_{m}\|_{H^{2}({\mathbb R}^{d})}\leq
\end{equation}
$$
c_{a}\|T_{m}(u_{0, m}(x)+v_{m}(x))\|_{H^{2}({\mathbb R}^{d})}
\Big\|\int_{{\mathbb R}^{d}}K_{m}(x-y)g_{m}(u_{0}(y)+v(y))dy\Big\|_{H^{2}({\mathbb R}^{d})}.
$$
Consider the right side of (\ref{auxub}). Clearly,
\begin{equation}
\label{tov}  
\|T_{m}(u_{0, m}(x)+v_{m}(x))\|_{H^{2}({\mathbb R}^{d})}\leq \|T_{m}\|
(\|u_{0}(x)\|_{H^{2}({\mathbb R}^{d}, {\mathbb R}^{N})}+1). 
\end{equation}
Let us recall inequality (\ref{y}). Hence,
$$
\Big\|\int_{{\mathbb R}^{d}}K_{m}(x-y)g_{m}(u_{0}(y)+v(y))dy\Big\|_{L^{2}({\mathbb R}^{d})}
\leq
$$
\begin{equation}
\label{kgvl2ub}
\|K_{m}\|_{L^{1}({\mathbb R}^{d})}\|g_{m}(u_{0}(x)+v(x))\|_{L^{2}({\mathbb R}^{d})}.
\end{equation}
Evidently, estimate (\ref{yd}) yields
$$
\Big\|\Delta_{x} \int_{{\mathbb R}^{d}}K_{m}(x-y)g_{m}(u_{0}(y)+v(y))dy\Big\|
_{L^{2}({\mathbb R}^{d})}\leq
$$
\begin{equation}
\label{kgvl2ubd}
\|\Delta K_{m} \|_{L^{1}({\mathbb R}^{d})}\|g_{m}(u_{0}(x)+v(x))\|_
{L^{2}({\mathbb R}^{d})}.
\end{equation}
Formulas (\ref{kgvl2ub}) and (\ref{kgvl2ubd}) give us that
$$
\Big\|\int_{{\mathbb R}^{d}}K_{m}(x-y)g_{m}(u_{0}(y)+v(y))dy\Big\|_{H^{2}({\mathbb R}^{d})}
\leq
$$
\begin{equation}
\label{kgvl2ubc}
\|K_{m}\|_{{\tilde W}^{2, 1}({\mathbb R}^{d})}\|g_{m}(u_{0}(x)+v(x))\|_{L^{2}({\mathbb R}^{d})}.
\end{equation}
Obviously, for $1\leq m\leq N$,
\begin{equation}
\label{gi}
g_{m}(u_{0}(x)+v(x))=\int_{0}^{1}
\nabla g_{m}(t(u_{0}(x)+v(x))).(u_{0}(x)+v(x))dt,
\end{equation}
where the dot denotes the scalar product of two vectors in ${\mathbb R}^{N}$.

Let us use bound (\ref{e}). Thus, for $v(x)\in B_{\rho}$ we derive
\begin{equation}
\label{u0pv}
|u_{0}+v|_{{\mathbb R}^{N}}\leq c_{e}
(\|u_{0}\|_{H^{2}({\mathbb R}^{d}, {\mathbb R}^{N})}+1). 
\end{equation}
Hence,
$$
|g_{m}(u_{0}(x)+v(x))|\leq
$$
\begin{equation}
\label{G}
\hbox{sup}_{z\in I}|\nabla g_{m}(z)|_{{\mathbb R}^{N}}
|u_{0}(x)+v(x)|_{{\mathbb R}^{N}}\leq
M|u_{0}(x)+v(x)|_{{\mathbb R}^{N}},
\end{equation}
where the ball $I$ is defined in (\ref{i}). Then
\begin{equation}
\label{Gn}
\|g_{m}(u_{0}(x)+v(x))\|_{L^{2}({\mathbb R}^{d})}
\leq M (\|u_{0}\|_{H^{2}({\mathbb R}^{d}, {\mathbb R}^{N})}+1).  
\end{equation}
By virtue of the estimates above we obtain that for $1\leq m\leq N$
\begin{equation}
\label{ur}  
\|u_{m}(x)\|_{H^{2}({\mathbb R}^{d})}\leq c_{a}M
(\|u_{0}\|_{H^{2}({\mathbb R}^{d}, {\mathbb R}^{N})}+1)^{2}
\|T_{m}\|\|K_{m}\|_{{\tilde W}^{2, 1}({\mathbb R}^{d})},
\end{equation}
so that
\begin{equation}
\label{urvf}  
\|u(x)\|_{H^{2}({\mathbb R}^{d}, {\mathbb R}^{N})}\leq c_{a}M
(\|u_{0}\|_{H^{2}({\mathbb R}^{d}, {\mathbb R}^{N})}+1)^{2}Q.
\end{equation}
Let us recall inequality (\ref{rub}). Hence,
$\|u(x)\|_{H^{2}({\mathbb R}^{d}, {\mathbb R}^{N})}\leq \rho$.
This means that the vector-function $u(x)$, which is uniquely determined by (\ref{aux}) 
belongs to $B_{\rho}$ as well.
Therefore, system (\ref{aux}) defines a map
$t_{g}: B_{\rho}\to B_{\rho}$ under the given conditions.

The goal is to establish that under the stated assumptions this map is a strict
contraction. We choose
arbitrarily $v_{1,2}(x)\in B_{\rho}$. According to the argument above,
$u_{1,2}:=t_{g}v_{1,2}\in B_{\rho}$ as well. By means of (\ref{aux}), for
$1\leq m\leq N, \ d=2,3$ we have
$$
u_{1, m}(x)=
$$
\begin{equation}
\label{aux1}
[T_{m}(u_{0, m}(x)+v_{1, m}(x))]\int_{{\mathbb R}^{d}}
K_{m}(x-y)g_{m}(u_{0}(y)+v_{1}(y))dy,
\end{equation}
$$
u_{2, m}(x)=
$$
\begin{equation}
\label{aux2}
[T_{m}(u_{0, m}(x)+v_{2, m}(x))]\int_{{\mathbb R}^{d}}
K_{m}(x-y)g_{m}(u_{0}(y)+v_{2}(y))dy.
\end{equation}
By virtue of (\ref{aux1}) and (\ref{aux2}),
\begin{equation}
\label{u1mu2} 
u_{1, m}(x)-u_{2, m}(x)=
\end{equation}
$$
[T_{m}v_{1, m}(x)-T_{m}v_{2, m}(x)]\int_{{\mathbb R}^{d}}K_{m}(x-y)
g_{m}(u_{0}(y)+v_{1}(y))dy+
$$
$$
[T_{m}(u_{0, m}(x)+v_{2, m}(x))]\times
$$
$$
\int_{{\mathbb R}^{d}}K_{m}(x-y)
[g_{m}(u_{0}(y)+v_{1}(y))-g_{m}(u_{0}(y)+v_{2}(y))]dy.
$$
From (\ref{u1mu2}) via  (\ref{alg}) we easily derive that
\begin{equation}
\label{u1mu2hin}
\|u_{1, m}(x)-u_{2, m}(x)\|_{H^{2}({\mathbb R}^{d})}\leq c_{a}\|T_{m}v_{1, m}(x)-
T_{m}v_{2, m}(x)\|_{H^{2}({\mathbb R}^{d})}\times
\end{equation}
$$  
\Big\|\int_{{\mathbb R}^{d}}K_{m}(x-y)g_{m}(u_{0}(y)+v_{1}(y))dy\Big\|_
{H^{2}({\mathbb R}^{d})}+
$$
$$
c_{a}\|T_{m}(u_{0, m}(x)+v_{2, m}(x))\|_{H^{2}({\mathbb R}^{d})}\times
$$
$$
\Big\|\int_{{\mathbb R}^{d}}K_{m}(x-y)
[g_{m}(u_{0}(y)+v_{1}(y))-g_{m}(u_{0}(y)+v_{2}(y))]dy\Big\|_{H^{2}({\mathbb R}^{d})}.
$$
Let us obtain the upper bound on the right side of inequality (\ref{u1mu2hin}). 
We have
\begin{equation}
\label{tv1mtv2}  
\|T_{m}v_{1, m}(x)-T_{m}v_{2, m}(x)\|_{H^{2}({\mathbb R}^{d})}\leq \|T_{m}\|
\|v_{1}(x)-v_{2}(x)\|_{H^{2}({\mathbb R}^{d}, {\mathbb R}^{N})}.
\end{equation}
Clearly, by the same argument as above     
$$
\Big\|\int_{{\mathbb R}^{d}}K_{m}(x-y)g_{m}(u_{0}(y)+v_{1}(y))dy\Big\|_
{H^{2}({\mathbb R}^{d})}\leq
$$
\begin{equation}
\label{kgv1intl2h1}  
\|K_{m}\|_{{\tilde W}^{2, 1}({\mathbb R}^{d})}M(\|u_{0}\|_{H^{2}({\mathbb R}^{d}, {\mathbb R}^{N})}+1).
\end{equation}
Therefore, the first term in the right side of inequality (\ref{u1mu2hin}) can
be estimated from above by
\begin{equation}
\label{catv12km}
c_{a}\|T_{m}\|\|v_{1}(x)-v_{2}(x)\|_{H^{2}({\mathbb R}^{d}, {\mathbb R}^{N})}
\|K_{m}\|_{{\tilde W}^{2, 1}({\mathbb R}^{d})}M
(\|u_{0}\|_{H^{2}({\mathbb R}^{d}, {\mathbb R}^{N})}+1).
\end{equation}
We consider the second term in the right side of
(\ref{u1mu2hin}). Evidently, for $1\leq m\leq N$,
\begin{equation}
\label{tu0v2h1}
\|T_{m}(u_{0, m}(x)+v_{2, m}(x))\|_{H^{2}({\mathbb R}^{d})}\leq \|T_{m}\|
(\|u_{0}\|_{H^{2}({\mathbb R}^{d}, {\mathbb R}^{N})}+1).  
\end{equation}
Let us recall inequalities (\ref{y}) and (\ref{yd}) to obtain that
$$
\Big\|\int_{{\mathbb R}^{d}}K_{m}(x-y)[g_{m}(u_{0}(y)+v_{1}(y))-
g_{m}(u_{0}(y)+v_{2}(y))]dy\Big\|_{H^{2}({\mathbb R}^{d})}\leq
$$
\begin{equation}
\label{kgu0v12l2h1}
\|K_{m}\|_{{\tilde W}^{2, 1}({\mathbb R}^{d})}\|
g_{m}(u_{0}(x)+v_{1}(x))-g_{m}(u_{0}(x)+v_{2}(x))\|_{L^{2}({\mathbb R}^{d})}.      
\end{equation}
Note that for $1\leq m\leq N$
$$
g_{m}(u_{0}(x)+v_{1}(x))-g_{m}(u_{0}(x)+v_{2}(x))=
$$
\begin{equation}
\label{gu0v12iz}
\int_{0}^{1}\nabla
g_{m}(u_{0}(x)+tv_{1}(x)+(1-t)v_{2}(x)).[v_{1}(x)-v_{2}(x)]dt.
\end{equation}  
Identity (\ref{gu0v12iz}) yields that
\begin{equation}
\label{gu0v12m}
|g_{m}(u_{0}(x)+v_{1}(x))-g_{m}(u_{0}(x)+v_{2}(x))|\leq
M|v_{1}(x)-v_{2}(x)|_{{\mathbb R}^{N}}.
\end{equation}
Then the norm
$$
\|g_{m}(u_{0}(x)+v_{1}(x))-g_{m}(u_{0}(x)+v_{2}(x))\|_{L^{2}({\mathbb R}^{d})}\leq
$$
\begin{equation}
\label{gu0v12n}
M\|v_{1}(x)-v_{2}(x)\|_{H^{2}({\mathbb R}^{d}, {\mathbb R}^{N})}.
\end{equation}
This means that the second term in the right side of  (\ref{u1mu2hin}) can
be also bounded from above by quantity (\ref{catv12km}). Therefore,
$$
\|u_{1}(x)-u_{2}(x)\|_{H^{2}({\mathbb R}^{d}, {\mathbb R}^{N})}\leq
$$
\begin{equation}
\label{u12v12h1}  
2c_{a}QM(\|u_{0}\|_{H^{2}({\mathbb R}^{d}, {\mathbb R}^{N})}+1)
\|v_{1}(x)-v_{2}(x)\|_{H^{2}({\mathbb R}^{d}, {\mathbb R}^{N})}.
\end{equation}
By means of inequality (\ref{u12v12h1}) and definition (\ref{sig}) we arrive at
\begin{equation}
\label{tgv12}
\|t_{g}v_{1}(x)-t_{g}v_{2}(x)\|_{H^{2}({\mathbb R}^{d}, {\mathbb R}^{N})}\leq \sigma 
\|v_{1}(x)-v_{2}(x)\|_{H^{2}({\mathbb R}^{d}, {\mathbb R}^{N})}.
\end{equation}
By virtue of condition (\ref{rub}) it can be trivially checked that the constant in the right
side of the estimate above
\begin{equation}
\label{sigm}  
\sigma<1.
\end{equation}
Therefore, the map $t_{g}: B_{\rho}\to B_{\rho}$ defined
by the system of equations (\ref{aux}) is a strict contraction under the given conditions.
Its unique fixed point $u_{p}(x)$ is the only solution of system 
(\ref{pert}) in the ball $B_{\rho}$. 
The resulting $u(x)$ given by formula (\ref{r}) solves our problem (\ref{p}).
\hspace{5cm} $\Box$

\bigskip

Let us conclude the work by establishing the validity of the second main
statement.

\bigskip


\setcounter{equation}{0}

\section{The continuity of the resulting solution with respect to
the vector function $g$}

\noindent
{\it Proof of Theorem 1.4.} Obviously, under the stated assumptions
\begin{equation}
\label{up1up2}  
u_{p,1}=t_{g_{1}}u_{p,1}, \quad u_{p,2}=t_{g_{2}}u_{p,2}.
\end{equation}
Then
\begin{equation}
\label{up1mup2}  
u_{p,1}-u_{p,2}=t_{g_{1}}u_{p,1}-t_{g_{1}}u_{p,2}+t_{g_{1}}u_{p,2}-
t_{g_{2}}u_{p,2},
\end{equation}
so that
$$
\|u_{p,1}-u_{p,2}\|_{H^{2}({\mathbb R}^{d}, {\mathbb R}^{N})}\leq
$$
\begin{equation}
\label{up1mup2n}  
\|t_{g_{1}}u_{p,1}-t_{g_{2}}u_{p,2}\|_{H^{2}({\mathbb R}^{d}, {\mathbb R}^{N})}+
\|t_{g_{1}}u_{p,2}-t_{g_{2}}u_{p,2}\|_{H^{2}({\mathbb R}^{d}, {\mathbb R}^{N})}.
\end{equation}
Let us recall estimate (\ref{tgv12}). Thus,
\begin{equation}
\label{tg1up1up2}  
\|t_{g_{1}}u_{p,1}-t_{g_{1}}u_{p,2}\|_{H^{2}({\mathbb R}^{d}, {\mathbb R}^{N})}\leq 
\sigma\|u_{p,1}-u_{p,2}\|_{H^{2}({\mathbb R}^{d}, {\mathbb R}^{N})}
\end{equation}
with $\sigma$ introduced in (\ref{sig}), such that bound (\ref{sigm}) holds.
Hence, we obtain
\begin{equation}
\label{sigma}
(1-\sigma)\|u_{p,1}-u_{p,2}\|_{H^{2}({\mathbb R}^{d}, {\mathbb R}^{N})}\leq
\|t_{g_{1}}u_{p,2}-t_{g_{2}}u_{p,2}\|_{H^{2}({\mathbb R}^{d}, {\mathbb R}^{N})}.
\end{equation}
Let us define
$\gamma(x):=t_{g_{1}}u_{p,2}$, so that for $1\leq m\leq N$
$$
\gamma_{m}(x)=
$$
\begin{equation}
\label{12}
[T_{m}(u_{0,m}(x)+u_{p,2,m}(x))]\int_{{\mathbb R}^{d}}
K_{m}(x-y)g_{1,m}(u_{0}(y)+u_{p,2}(y))dy,
\end{equation}
$$
u_{p,2,m}(x)=
$$
\begin{equation}
\label{22}
[T_{m}(u_{0,m}(x)+u_{p,2,m}(x))]\int_{{\mathbb R}^{d}}
K_{m}(x-y)g_{2,m}(u_{0}(y)+u_{p,2}(y))dy.
\end{equation}
By virtue of formulas (\ref{12}) and (\ref{22}), we have
$$
\gamma_{m}(x)-u_{p,2,m}(x)=[T_{m}(u_{0,m}(x)+u_{p,2,m}(x))]\times
$$
\begin{equation}
\label{rup2}  
\int_{{\mathbb R}^{d}}K_{m}(x-y)
[g_{1,m}(u_{0}(y)+u_{p,2}(y))-g_{2,m}(u_{0}(y)+u_{p,2}(y))]dy.
\end{equation}
Using (\ref{alg}), we derive
$$
\|\gamma_{m}(x)-u_{p,2,m}(x)\|_{H^{2}({\mathbb R}^{d})}\leq c_{a}
\|T_{m}(u_{0,m}(x)+u_{p,2,m}(x))\|_{H^{2}({\mathbb R}^{d})}\times
$$
\begin{equation}
\label{rup2n}
\Big\|\int_{{\mathbb R}^{d}}K_{m}(x-y)
[g_{1,m}(u_{0}(y)+u_{p,2}(y))-g_{2,m}(u_{0}(y)+u_{p,2}(y))]dy\Big\|_{H^{2}({\mathbb R}^{d})}.
\end{equation}
Clearly, the inequality
\begin{equation}
\label{tu0up2n}  
\|T_{m}(u_{0,m}(x)+u_{p,2,m}(x))\|_{H^{2}({\mathbb R}^{d})}\leq \|T_{m}\|
(\|u_{0}\|_{H^{2}({\mathbb R}^{d}, {\mathbb R}^{N})}+1)
\end{equation}
is valid.
Let us recall estimates (\ref{y}) and (\ref{yd}). This yields
$$
\Big\|\int_{{\mathbb R}^{d}}K_{m}(x-y)
[g_{1,m}(u_{0}(y)+u_{p,2}(y))-g_{2,m}(u_{0}(y)+u_{p,2}(y))]dy\Big\|_{H^{2}({\mathbb R}^{d})}
\leq
$$
\begin{equation}
\label{kg1g2u0up2nh}  
\|K_{m}\|_{{\tilde W}^{2, 1}({\mathbb R}^{d})}\|g_{1,m}(u_{0}(x)+u_{p,2}(x))-
g_{2,m}(u_{0}(x)+u_{p,2}(x))\|_{L^{2}({\mathbb R}^{d})}.
\end{equation}
Obviously, for $1\leq m\leq N$,
\begin{equation}
\label{g12u0up2}  
[g_{1,m}-g_{2,m}](u_{0}(x)+u_{p,2}(x))=
\end{equation}
$$
\int_{0}^{1}\nabla[g_{1,m}-g_{2,m}](t(u_{0}(x)+u_{p,2}(x))).[u_{0}(x)+u_{p,2}(x)]dt.
$$
By means of equality (\ref{g12u0up2}) we obtain
\begin{equation}
\label{g12u0up2m}  
|[g_{1,m}-g_{2,m}](u_{0}(x)+u_{p,2}(x))|\leq
\|g_{1,m}-g_{2,m}\|_{C^{1}(I)}|u_{0}(x)+u_{p,2}(x)|_{{\mathbb R}^{N}}.
\end{equation}
Let us derive the upper bound for the norm as
$$
\|[g_{1,m}-g_{2,m}](u_{0}(x)+u_{p,2}(x))\|_{L^{2}({\mathbb R}^{d})}\leq
$$
\begin{equation}
\label{g12u0up2mn}
\|g_{1}-g_{2}\|_{C^{1}(I, {\mathbb R}^{N})}(\|u_{0}\|_{H^{2}({\mathbb R}^{d}, {\mathbb R}^{N})}+1).
\end{equation}
Using inequalities (\ref{rup2n}), (\ref{tu0up2n}), 
(\ref{kg1g2u0up2nh}), (\ref{g12u0up2mn}),  we arrive at
$$
\|\gamma(x)-u_{p,2}(x)\|_{H^{2}({\mathbb R}^{d}, {\mathbb R}^{N})}\leq
$$
\begin{equation}
\label{rup2h1r}
c_{a}Q
(\|u_{0}\|_{H^{2}({\mathbb R}^{d}, {\mathbb R}^{N})}+1)^{2}\|g_{1}(z)-g_{2}(z)\|_
{C^{1}(I, {\mathbb R}^{N})}.
\end{equation}
According to formulas (\ref{sigma}) and (\ref{rup2h1r}), we obtain that
$$
\|u_{p,1}(x)-u_{p,2}(x)\|_{H^{2}({\mathbb R}^{d}, {\mathbb R}^{N})}\leq
$$
\begin{equation}
\label{p1p2}
\frac{c_{a}}{1-\sigma}
(\|u_{0}\|_{H^{2}({\mathbb R}^{d}, {\mathbb R}^{N})}+1)^{2}Q
\|g_{1}(z)-g_{2}(z)\|_{C^{1}(I, {\mathbb R}^{N})}.
\end{equation}
By virtue of (\ref{u12}) along with estimate (\ref{p1p2}) and definition
(\ref{sig}), inequality (\ref{cont}) holds.  \hspace{9.5cm} $\Box$

\bigskip


\section{Acknowledgement} The work was partially supported by the
NSERC Discovery grant.

\bigskip


\end{document}